\documentstyle{article}
\input amssym.def
\input amssym
\input epsf.tex

\begin{document}
\newtheorem{thm}{Theorem}[section]
\newtheorem{defin}[thm]{Definition}
\newtheorem{lemma}[thm]{Lemma}
\newtheorem{propo}[thm]{Proposition}
\newtheorem{cor}[thm]{Corollary}
\newtheorem{conj}[thm]{Conjecture}

\centerline{\LARGE \bf Generalized Barrett-Crane Vertices}

\centerline{\LARGE \bf and Invariants of Embedded Graphs}
\bigskip

\centerline{\parbox{56mm}{David N. Yetter \\ Department of
Mathematics \\ Kansas State University \\ Manhattan, KS 66506}
\footnote{Supported by NSF Grant \# DMS-9504423}}
\bigskip

\section {Introduction}

In \cite{BC} Barrett and Crane introduce a modification of the generalized
Crane-Yetter state-sum (cf. \cite{CKY}) based on the category of 
representations of $Spin(4) \cong SU(2)\times SU(2)$, which provides
a four-dimensional analogue of Regge and Ponzano's \cite{RP} spin-network 
formulation
of three-dimensional gravity.

The key to the modification of the Crane-Yetter state-sum is the use
of a different intertwiner between the ``inbound'' and ``outbound''
tensor products of objects assigned to faces of the tetrahedron, thereby
imposing a ``quantum analogue'' of the condition that the sum of the
simple bivectors represented by two faces with a common edge is itself
simple. The 
purposes of this paper are

\begin{enumerate}
\item to given an explicit formula for the Barrett-Crane
intertwiners and their $q$-analogues by using 
the Kauffman-Lins \cite{KL} ``Temperley-Lieb'' 
recoupling theory for (quantized) $SU(2)$ in each factor, 
\item to show
that in the case where the $SU(2)$ factors of $Spin(4)$ are deformed
with reciprocal deformation parameters the operators have the symmetry
properties of a topological 4-vertex when represented in the usual 
diagrammatic calculus for braided tensor categories (cf. \cite{JS}, \cite{FY}
\cite{Y.rep}), 
\item to generalize the family of operators to families which exhibit
the symmetry properties of topological $n$-vertices when represented in the
diagrammatic calculus, and
\item to consider the natural family of invariants of topological graphs
embedded in ${\Bbb R}^3$ or $S^3$ which arise by
applying the functorial construction
 first proposed in \cite{FY} and \cite{Y.graph}.    
\end{enumerate}

Although the present paper will be concerned with the construction of
``algebraic $n$-vertices'' and their applications to topological graph theory,
the operators constructed herein may well have applications to the
physically-motivated Barrett-Crane state-sum, in particular allowing 
computations of state spaces and transition amplitudes from cell-decompositions
coarser than triangulations.

\newpage
\section{$Spin(4)$ Recoupling Theory via Temperley-Lieb Recoupling}

The factorization $Spin(4) \cong SU(2) \times SU(2)$ (or rather
$\frak{so}_4 \cong \frak{sl}_2 \oplus \frak{sl}_2$) has an easy, immediate
consequence for the structure of the associated QUEA:  the standard
1-parameter deformation $U_q(\frak{sl}_2)$ gives rise to a 2-parameter
deformation $U_{q,r}(\frak{so}_4) \cong U_q(\frak{sl}_2)\otimes 
U_r(\frak{sl}_2)$, where $q$ is the deformation parameter
for the first $\frak{sl}_2$, and $r$ for the second.

	We will be interested in 1-parameter deformation which arises by
specializing to $r = q^{-1}$, which we will denote 
$U_q^{\rm bal}(\frak{so}_4)$. The irreducible representations of
$U_q^{\rm bal}(\frak{so}_4)$ are each a tensor product of an irreducible 
representation of $U_q(\frak{sl}_2)$ with an irreducible
representation of $U_{q^{-1}}(\frak{sl}_2)$, and thus are indexed by
pairs of non-negative half-integers (spins), or more coveniently, since
we will be using Kauffman-Lins style recombination by integers (twice
spin = dimension $- 1$ = number of strands $q$-symmetrized).

	Throughout the following we adopt the notation of Kauffman-Lins
\cite{KL} for all aspect of $U_q(\frak{sl}_2)$ recoupling theory.  It should
be observed that it will be unnecessary to have
different notations for $U_q(\frak{sl}_2)$ and $U_{q^{-1}}(\frak{sl}_2)$
recoupling constants, since the $q$-symmetrizer, three-vertex, quantum
dimensions, $\theta(a,b,c)$, and $q$-6j symbols are all unchanged when
$q$ is replaced with $q^{-1}$.

	We can then let a recombination network labelled with pairs of
twice-spins represent the tensor product of the corresponding network
labelled with first twice-spins with that labelled with second twice-spins, 
preceeded
and followed by the necessary coherence maps from the underlying 
category of vector-spaces.

	For example

\begin{figure}[h]
\begin{center}
\setlength{\unitlength}{0.00083300in}%
\begingroup\makeatletter\ifx\SetFigFont\undefined
% extract first six characters in \fmtname
\def\x#1#2#3#4#5#6#7\relax{\def\x{#1#2#3#4#5#6}}%
\expandafter\x\fmtname xxxxxx\relax \def\y{splain}%
\ifx\x\y   % LaTeX or SliTeX?
\gdef\SetFigFont#1#2#3{%
  \ifnum #1<17\tiny\else \ifnum #1<20\small\else
  \ifnum #1<24\normalsize\else \ifnum #1<29\large\else
  \ifnum #1<34\Large\else \ifnum #1<41\LARGE\else
     \huge\fi\fi\fi\fi\fi\fi
  \csname #3\endcsname}%
\else
\gdef\SetFigFont#1#2#3{\begingroup
  \count@#1\relax \ifnum 25<\count@\count@25\fi
  \def\x{\endgroup\@setsize\SetFigFont{#2pt}}%
  \expandafter\x
    \csname \romannumeral\the\count@ pt\expandafter\endcsname
    \csname @\romannumeral\the\count@ pt\endcsname
  \csname #3\endcsname}%
\fi
\fi\endgroup
\begin{picture}(1512,1248)(901,-1072)
\thicklines
\put(1801,-436){\circle*{150}}
\put(1201,164){\line( 1,-1){600}}
\put(1801,-436){\line( 1, 1){600}}
\put(1801,-436){\line( 0,-1){600}}
\put(901,-136){\makebox(0,0)[lb]{\smash{\SetFigFont{10}{14.4}{rm}$i,j$}}}
\put(2401,-211){\makebox(0,0)[lb]{\smash{\SetFigFont{10}{14.4}{rm}$k,l$}}}
\put(1951,-1036){\makebox(0,0)[lb]{\smash{\SetFigFont{10}{14.4}{rm}$m,n$}}}
\end{picture}

\end{center}	
\end{figure}

\noindent represents the map  

\[ (V_i\otimes W_j)\otimes (V_k\otimes W_l)\stackrel{c}{\rightarrow} 
(V_i\otimes V_k) \otimes (W_j\otimes W_l)
\stackrel{v_{ikm}\otimes w_{jln}}{\longrightarrow} V_m \otimes W_n \]

\noindent where $c$ is the middle-four-interchange map composed of 
associators and symmetrizers in ${\Bbb C}-v.s.$, and $v_{ikm}$ (resp.
$w_{jln}$) is the map named by the Kauffman-Lins three-vertex with the
chosen labels.

	By the same reasoning, it follows that the quantum dimension
of $(i,j)$ is $\Delta_i \Delta_j$, and similarly that the evaluation of
a closed diagram labelled with pairs of twice-spins is the product of 
the evaluations of the corresponding diagrams obtained by selecting
all of the first (resp. second) twice-spins as labels.  (In the case
where there are braidings, one must use the correct deformation
parameter in each case).

	A crucial ingredient in the Barrett-Crane state sum and in 
what follows is the notion of a 
``balanced'' irreducible, that is one of the form $(j,j)$ for some
$j$.  Several facts stand out about balanced irreducibles in recombination
diagrams:

First, it should be observed that a curl on a strand with a balanced 
label is simply the identity operator:  the constants contributed by
each tensorand cancel (because we are working with 
$U_q^{\rm bal}(\frak{so}_4)$) .

Second, and for the same reason, we have

\begin{figure}[h]
\begin{center}
\setlength{\unitlength}{0.00083300in}%
\begingroup\makeatletter\ifx\SetFigFont\undefined
% extract first six characters in \fmtname
\def\x#1#2#3#4#5#6#7\relax{\def\x{#1#2#3#4#5#6}}%
\expandafter\x\fmtname xxxxxx\relax \def\y{splain}%
\ifx\x\y   % LaTeX or SliTeX?
\gdef\SetFigFont#1#2#3{%
  \ifnum #1<17\tiny\else \ifnum #1<20\small\else
  \ifnum #1<24\normalsize\else \ifnum #1<29\large\else
  \ifnum #1<34\Large\else \ifnum #1<41\LARGE\else
     \huge\fi\fi\fi\fi\fi\fi
  \csname #3\endcsname}%
\else
\gdef\SetFigFont#1#2#3{\begingroup
  \count@#1\relax \ifnum 25<\count@\count@25\fi
  \def\x{\endgroup\@setsize\SetFigFont{#2pt}}%
  \expandafter\x
    \csname \romannumeral\the\count@ pt\expandafter\endcsname
    \csname @\romannumeral\the\count@ pt\endcsname
  \csname #3\endcsname}%
\fi
\fi\endgroup
\begin{picture}(3912,2145)(901,-2269)
\thicklines
\put(1801,-1336){\circle*{150}}
\put(4276,-1261){\circle*{150}}
\put(1201,-136){\line( 1,-1){900}}
\put(2101,-1036){\line(-1,-1){300}}
\put(1801,-1336){\line(-1, 1){300}}
\put(1501,-1036){\line( 1, 1){225}}
\put(1876,-661){\line( 1, 1){525}}
\put(1801,-1336){\line( 0,-1){900}}
\put(4276,-1261){\line( 0,-1){900}}
\put(3601,-136){\line( 3,-5){675}}
\put(4276,-1261){\line( 1, 2){525}}
\put(2851,-1261){\makebox(0,0)[lb]{\smash{\SetFigFont{10}{14.4}{rm}=}}}
\put(901,-361){\makebox(0,0)[lb]{\smash{\SetFigFont{10}{14.4}{rm}$i,i$}}}
\put(2476,-436){\makebox(0,0)[lb]{\smash{\SetFigFont{10}{14.4}{rm}$j,j$}}}
\put(2026,-2236){\makebox(0,0)[lb]{\smash{\SetFigFont{10}{14.4}{rm}$k,k$}}}
\put(3376,-361){\makebox(0,0)[lb]{\smash{\SetFigFont{10}{14.4}{rm}$i,i$}}}
\put(4801,-436){\makebox(0,0)[lb]{\smash{\SetFigFont{10}{14.4}{rm}$j,j$}}}
\put(4501,-2161){\makebox(0,0)[lb]{\smash{\SetFigFont{10}{14.4}{rm}$k,k$}}}
\end{picture}

\end{center}
\end{figure}

	Thus we see that for balanced irreducibles, the curl is trivial,
while the (family of) 3-vertex admits as symmetries all permuations, 
not just cyclic permutations as in Kauffman-Lins (or for non-balanced
labels in the present setting).

	We are now in a position to introduce the formula in Kauffman-Lins
notation for the Barrett-Crane 4-vertex.

\newpage
\section{Barrett-Crane Vertices, Properties and Generalizations}

\begin{defin}
The {\em Barrett-Crane 4-vertex} $v^{i,j}_{k,l}$ is the map
from $V_{(i,i)}\otimes V_{(j,j)}$ to $V_{(k,k)}\otimes V_{(l,l)}$ given
by the sum of recombination networks
\end{defin}

\begin{figure}[h]
\begin{center}

\setlength{\unitlength}{0.00083300in}%
\begingroup\makeatletter\ifx\SetFigFont\undefined
% extract first six characters in \fmtname
\def\x#1#2#3#4#5#6#7\relax{\def\x{#1#2#3#4#5#6}}%
\expandafter\x\fmtname xxxxxx\relax \def\y{splain}%
\ifx\x\y   % LaTeX or SliTeX?
\gdef\SetFigFont#1#2#3{%
  \ifnum #1<17\tiny\else \ifnum #1<20\small\else
  \ifnum #1<24\normalsize\else \ifnum #1<29\large\else
  \ifnum #1<34\Large\else \ifnum #1<41\LARGE\else
     \huge\fi\fi\fi\fi\fi\fi
  \csname #3\endcsname}%
\else
\gdef\SetFigFont#1#2#3{\begingroup
  \count@#1\relax \ifnum 25<\count@\count@25\fi
  \def\x{\endgroup\@setsize\SetFigFont{#2pt}}%
  \expandafter\x
    \csname \romannumeral\the\count@ pt\expandafter\endcsname
    \csname @\romannumeral\the\count@ pt\endcsname
  \csname #3\endcsname}%
\fi
\fi\endgroup
\begin{picture}(4362,2211)(1789,-2410)
\thicklines
\put(5401,-811){\circle*{150}}
\put(5401,-1786){\circle*{150}}
\put(2401,-811){\line(-1, 0){600}}
\put(1801,-811){\line( 1,-1){375}}
\put(2176,-1186){\line(-1,-1){375}}
\put(1801,-1561){\line( 1, 0){675}}
\put(2401,-811){\line( 1, 0){ 75}}
\put(2851,-1186){\line( 1, 0){1275}}
\put(4801,-211){\line( 1,-1){600}}
\put(5401,-811){\line( 1, 1){600}}
\put(4801,-2386){\line( 1, 1){600}}
\put(5401,-1786){\line( 1,-1){600}}
\put(5401,-811){\line( 0,-1){1050}}
\put(2026,-1861){\makebox(0,0)[lb]{\smash{\SetFigFont{10}{14.4}{rm}$n$}}}
\put(3301,-961){\makebox(0,0)[lb]{\smash{\SetFigFont{10}{14.4}{rm}$\Delta_n$}}}
\put(2926,-1561){\makebox(0,0)[lb]
{\smash{\SetFigFont{10}{14.4}{rm}$\theta(m,l,n) \theta(j,k,n)$}}}
\put(4501,-436){\makebox(0,0)[lb]{\smash{\SetFigFont{10}{14.4}{rm}$j,j$}}}
\put(6001,-511){\makebox(0,0)[lb]{\smash{\SetFigFont{10}{14.4}{rm}$k,k$}}}
\put(5551,-1411){\makebox(0,0)[lb]{\smash{\SetFigFont{10}{14.4}{rm}$n,n$}}}
\put(4501,-2311){\makebox(0,0)[lb]{\smash{\SetFigFont{10}{14.4}{rm}$l,l$}}}
\put(6151,-2386){\makebox(0,0)[lb]{\smash{\SetFigFont{10}{14.4}{rm}$m,m$}}}
\end{picture}

\end{center}
\end{figure}

From the point of view of this paper the crucial properties of these map are
their symmetry properties.  It is immediate from the observations above
and the cyclic symmetry properties of 3-vertices (which follow from
the corresponding properties for Kauffman-Lins 3-vertices) that

\[ \sigma_{(j,j),(i,i)} (v^{i,j}_{k,l}) = v^{j,i}_{k,l} \]

\[  v^{i,j}_{k,l}(\sigma_{(k,k),(l,l)}) = v^{i,j}_{l,k} \]

\[ \epsilon_{(i,i)\otimes (j,j)}\otimes V_{(l,l)}\otimes V_{(k,k)}(
V_{(j,j)}\otimes V_{i,i}\otimes v^{i,j}_{k,l} \otimes 
V_{(l,l)}\otimes V_{(k,k)}(V_{(j,j)}\otimes V_{(i,i)}\otimes 
\eta_{(k,k)\otimes (l,l)})) \\ \vspace*{1cm} = v^{l,k}_{j,i} \]

\noindent and

\[ V_{(l,l)}\otimes V_{(k,k)}\otimes \epsilon_{(i,i)\otimes (j,j)}(
V_{(l,l)}\otimes V_{(k,k)}\otimes v^{i,j}_{k,l}\otimes V_{(j,j)}\otimes V_{i,i}
(\eta_{(k,k)\otimes (l,l)}\otimes V_{(j,j)}\otimes V_{(i,i)})) 
\\ \vspace*{1cm} = v^{l,k}_{j,i} \]

\noindent where $\eta_{(-)}$ and $\epsilon_{(-)}$ are the duality 
transformations in $Rep(U_q^{\rm bal}(\frak{so}_4))$ normalized so that
they are given by the ``cup'' and ``cap'' networks, and intervening
generalized associators have been omitted by Mac Lane's coherence theorem
\cite{CWM}.

What is remarkable is that the Barrett-Crane 4-vertex can also be expressed
in terms of the same formula with the network rotated by $\pi /2$, that is:

\begin{propo}\label{turning.prop}

\[ v^{i,j}_{k,l} = \epsilon_{(k,k)}\otimes V_{(i,i)}\otimes V_{(j,j)}
(V_{(k,k)}\otimes v^{k,i}_{l,j}\otimes V_{(j,j)}(V_{(k,k)}\otimes V_{(l,l)}
\otimes \eta_{(j,j)})) \]

\noindent where $\epsilon$ and $\eta$ are as above, and the intervening
associators have been omitted.
\end{propo}

\begin{figure}[hbt]
\begin{center}
\epsfxsize=4.75in \epsfbox{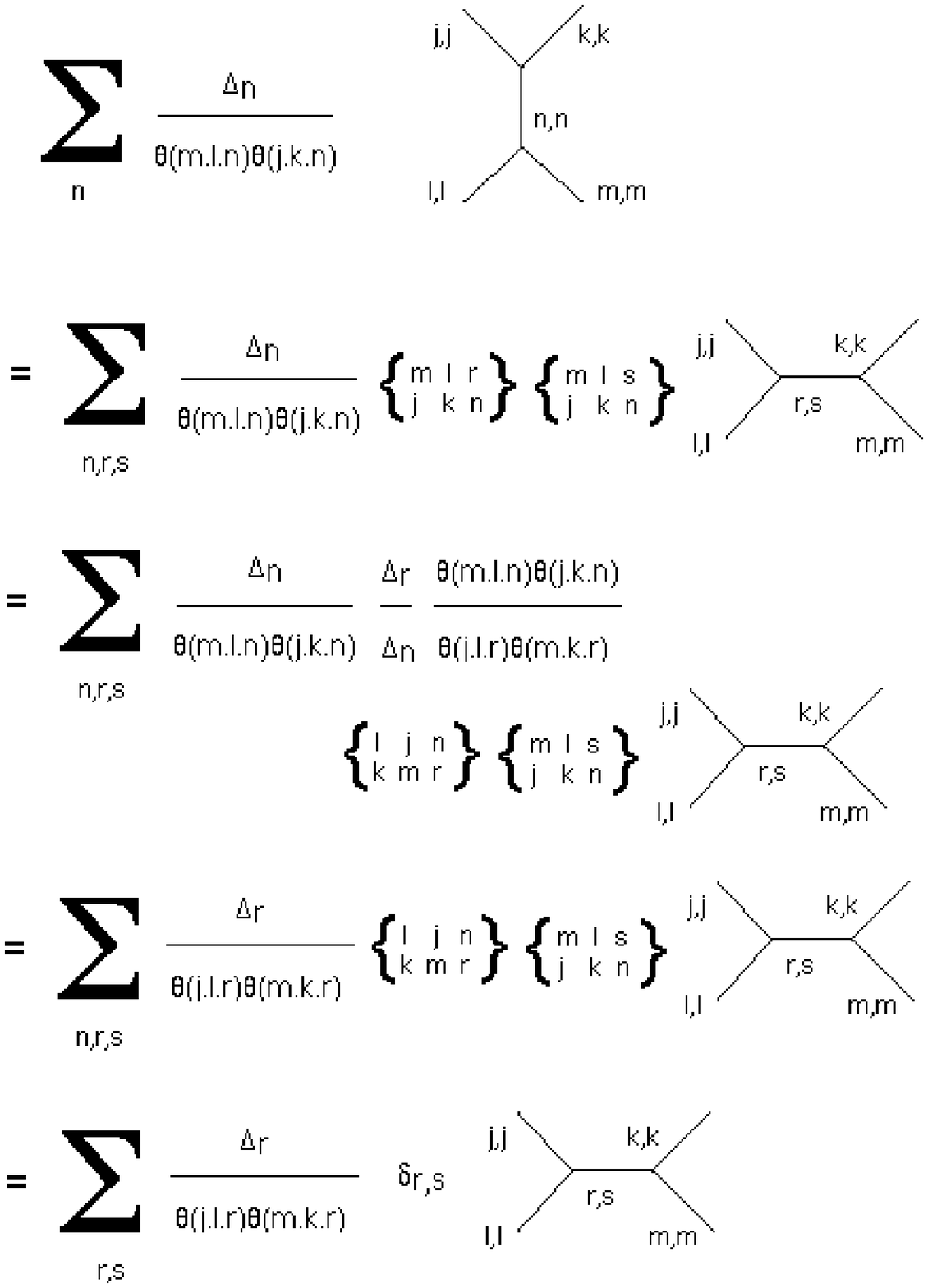}
\caption{\label{turn} Proof of Proposition \ref{turning.prop}}
\end{center}
\end{figure}

\noindent{\bf proof:} The calculation with linear combinations of 
recombination networks is given in Figure \ref{turn}, where $\delta_{r,s}$ is
the Kronecker delta.  The first equation
follows from the fact that the evaluation of a doubled planar network
(including all those in the Kauffman-Lins construction of a q-6j symbol)
is simply the product of the two corresponding $U_q({\frak sl}_2)$
networks, the second by the symmetry properties of the q-6j symbol,
the third is elementary, the fourth by the orthogonality properties of
the q-6j symbols.$\Box$

\smallskip
We will also refer to maps obtained from Barrett-Crane 4-vertices by
tensoring $V^{i,j}_{k,l}$ on the right (resp. right, left, left)
with $(j,j)$ (resp. $(l,l)$, $(i,i)$, $(k,k)$) 
and pre- (resp. post-, pre-, post-)composing 
with $V_{(i,i)}\otimes \eta_{(j,j)}$ (resp. 
$V_{(k,k)}\otimes \epsilon_{(l,l)}$,
$\eta_{(i,i)}\otimes V_{(j,j)}$,  $\epsilon_{(k,k}\otimes V_{(l,l)}$)
as Barrett-Crane 4-vertices, as well as those obtained by tensoring
$V^{i,j}_{k,l}$ on the right (resp. right, left, left) with $(j,i)$ (resp.
$(l,k)$, $(j,i)$, $(l,k)$) and pre- (resp. post-, pre- post-)composing
with $\eta_{(i,i)\otimes (j,j)}$ (resp. $\epsilon_{(l,l)\otimes (k,k)}$,
$\eta_{(j,j)\otimes (i,i)}$, $\epsilon_{(k,k)\otimes (l,l)}$).
Less formally, but more intelligibly, the maps just described are 
those obtained by using duality
to ``turn'' some of the inputs or outputs of the vertex ``down'' or ``up''--
the geometry of the spin-networks coinciding nicely with raising and lowering
of indices.

It thus follows from Proposition \ref{turning.prop} and the symmetry
properties already observed any map in 
$Rep(U_q^{\rm bal}(\frak{so}_4))$ between tensor products of   
balanced objects which admits an expression in terms of a connected
recombination network containing a single Barrett-Crane 4-vertex
is itself a Barrett-Crane 4-vertex in this more general sense.

Before generalizing to $n$-vertices, it is conveinent to change notation
slightly.  Let an unmarked
node where three edges labelled $(i,i)$, $(j,j)$ and $(k,k)$
meet denote $\frac{1}{\theta(i,j,k)}$ times the 3-vertex denoted with the heavy
dot:

\begin{figure}[h]
\setlength{\unitlength}{0.00083300in}%
\begingroup\makeatletter\ifx\SetFigFont\undefined
% extract first six characters in \fmtname
\def\x#1#2#3#4#5#6#7\relax{\def\x{#1#2#3#4#5#6}}%
\expandafter\x\fmtname xxxxxx\relax \def\y{splain}%
\ifx\x\y   % LaTeX or SliTeX?
\gdef\SetFigFont#1#2#3{%
  \ifnum #1<17\tiny\else \ifnum #1<20\small\else
  \ifnum #1<24\normalsize\else \ifnum #1<29\large\else
  \ifnum #1<34\Large\else \ifnum #1<41\LARGE\else
     \huge\fi\fi\fi\fi\fi\fi
  \csname #3\endcsname}%
\else
\gdef\SetFigFont#1#2#3{\begingroup
  \count@#1\relax \ifnum 25<\count@\count@25\fi
  \def\x{\endgroup\@setsize\SetFigFont{#2pt}}%
  \expandafter\x
    \csname \romannumeral\the\count@ pt\expandafter\endcsname
    \csname @\romannumeral\the\count@ pt\endcsname
  \csname #3\endcsname}%
\fi
\fi\endgroup
\begin{picture}(6225,1395)(976,-1594)
\thicklines
\put(1201,-211){\line( 1,-1){600}}
\put(1801,-811){\line( 1, 1){600}}
\put(1801,-811){\line( 0,-1){750}}
\put(976,-511){\makebox(0,0)[lb]{\smash{\SetFigFont{12}{14.4}{rm}$i,i$}}}
\put(2476,-511){\makebox(0,0)[lb]{\smash{\SetFigFont{12}{14.4}{rm}$j,j$}}}
\put(2026,-1561){\makebox(0,0)[lb]{\smash{\SetFigFont{12}{14.4}{rm}$k,k$}}}
\put(5926,-211){\line( 1,-1){600}}
\put(6526,-811){\line( 1, 1){600}}
\put(6526,-811){\line( 0,-1){750}}
\put(5701,-511){\makebox(0,0)[lb]{\smash{\SetFigFont{12}{14.4}{rm}$i,i$}}}
\put(7201,-511){\makebox(0,0)[lb]{\smash{\SetFigFont{12}{14.4}{rm}$j,j$}}}
\put(6751,-1561){\makebox(0,0)[lb]{\smash{\SetFigFont{12}{14.4}{rm}$k,k$}}}
\put(6526,-811){\circle*{150}}
\put(3676,-961){\makebox(0,0)[lb]{\smash{\SetFigFont{12}{14.4}{rm}$=$}}}
\put(4351,-1036){\makebox(0,0)[lb]
{\smash{\SetFigFont{12}{14.4}{rm}$\frac{1}{\theta(i,j,k)}$}}}
\end{picture}

\end{figure}

\noindent In this new notation the formula for a Barrett-Crane 4-vertex
becomes

\begin{figure}[h]
\begin{center}

\setlength{\unitlength}{0.00083300in}%
\begingroup\makeatletter\ifx\SetFigFont\undefined
% extract first six characters in \fmtname
\def\x#1#2#3#4#5#6#7\relax{\def\x{#1#2#3#4#5#6}}%
\expandafter\x\fmtname xxxxxx\relax \def\y{splain}%
\ifx\x\y   % LaTeX or SliTeX?
\gdef\SetFigFont#1#2#3{%
  \ifnum #1<17\tiny\else \ifnum #1<20\small\else
  \ifnum #1<24\normalsize\else \ifnum #1<29\large\else
  \ifnum #1<34\Large\else \ifnum #1<41\LARGE\else
     \huge\fi\fi\fi\fi\fi\fi
  \csname #3\endcsname}%
\else
\gdef\SetFigFont#1#2#3{\begingroup
  \count@#1\relax \ifnum 25<\count@\count@25\fi
  \def\x{\endgroup\@setsize\SetFigFont{#2pt}}%
  \expandafter\x
    \csname \romannumeral\the\count@ pt\expandafter\endcsname
    \csname @\romannumeral\the\count@ pt\endcsname
  \csname #3\endcsname}%
\fi
\fi\endgroup
\begin{picture}(4362,2211)(1789,-2410)
\thicklines
\put(2401,-811){\line(-1, 0){600}}
\put(1801,-811){\line( 1,-1){375}}
\put(2176,-1186){\line(-1,-1){375}}
\put(1801,-1561){\line( 1, 0){675}}
\put(2401,-811){\line( 1, 0){ 75}}
%\put(2851,-1186){\line( 1, 0){1275}}
\put(4801,-211){\line( 1,-1){600}}
\put(5401,-811){\line( 1, 1){600}}
\put(4801,-2386){\line( 1, 1){600}}
\put(5401,-1786){\line( 1,-1){600}}
\put(5401,-811){\line( 0,-1){1050}}
\put(2026,-1861){\makebox(0,0)[lb]{\smash{\SetFigFont{12}{14.4}{rm}$n$}}}
\put(3301,-961){\makebox(0,0)[lb]{\smash{\SetFigFont{12}{14.4}{rm}$\Delta_n$}}}
\put(4501,-436){\makebox(0,0)[lb]{\smash{\SetFigFont{12}{14.4}{rm}$j,j$}}}
\put(6001,-511){\makebox(0,0)[lb]{\smash{\SetFigFont{12}{14.4}{rm}$k,k$}}}
\put(5551,-1411){\makebox(0,0)[lb]{\smash{\SetFigFont{12}{14.4}{rm}$n,n$}}}
\put(4501,-2311){\makebox(0,0)[lb]{\smash{\SetFigFont{12}{14.4}{rm}$l,l$}}}
\put(6151,-2386){\makebox(0,0)[lb]{\smash{\SetFigFont{12}{14.4}{rm}$m,m$}}}
\end{picture}

\end{center}
\end{figure}  
\newpage

With this change of notation, it then becomes clear how to define $n$-vertices
for any $n \geq 3$.

For $n=3$ the 3-vertex is the ``new 3-vertex''

\begin{figure}[h]
\begin{center}
\setlength{\unitlength}{0.00083300in}%
\begingroup\makeatletter\ifx\SetFigFont\undefined
% extract first six characters in \fmtname
\def\x#1#2#3#4#5#6#7\relax{\def\x{#1#2#3#4#5#6}}%
\expandafter\x\fmtname xxxxxx\relax \def\y{splain}%
\ifx\x\y   % LaTeX or SliTeX?
\gdef\SetFigFont#1#2#3{%
  \ifnum #1<17\tiny\else \ifnum #1<20\small\else
  \ifnum #1<24\normalsize\else \ifnum #1<29\large\else
  \ifnum #1<34\Large\else \ifnum #1<41\LARGE\else
     \huge\fi\fi\fi\fi\fi\fi
  \csname #3\endcsname}%
\else
\gdef\SetFigFont#1#2#3{\begingroup
  \count@#1\relax \ifnum 25<\count@\count@25\fi
  \def\x{\endgroup\@setsize\SetFigFont{#2pt}}%
  \expandafter\x
    \csname \romannumeral\the\count@ pt\expandafter\endcsname
    \csname @\romannumeral\the\count@ pt\endcsname
  \csname #3\endcsname}%
\fi
\fi\endgroup
\begin{picture}(1512,1248)(901,-1072)
\thicklines
\put(1201,164){\line( 1,-1){600}}
\put(1801,-436){\line( 1, 1){600}}
\put(1801,-436){\line( 0,-1){600}}
\put(901,-136){\makebox(0,0)[lb]{\smash{\SetFigFont{12}{14.4}{rm}$i,i$}}}
\put(2401,-211){\makebox(0,0)[lb]{\smash{\SetFigFont{12}{14.4}{rm}$j,j$}}}
\put(1951,-1036){\makebox(0,0)[lb]{\smash{\SetFigFont{12}{14.4}{rm}$k,k$}}}
\end{picture}

\end{center}	
\end{figure}

\noindent for $n=4$ it is the Barrett-Crane 4-vertex already defined, while
for $n > 4$, an $n$-vertex with given
tensor product of $k$ balanced irreducibles
as source, and given tensor product of $n-k$ balanced irreducibles as target
can be described by a sum of recombination networks with 
underlying graph a fixed
tree with $n$ leaves (or better still $n$ ``external edges'' with a 
free end incident with a vertex) divided into a set of $k$ at the top,
and a set of $n-k$ at the bottom and all internal vertices trivalent,
with the external edges labelled with the
tensorands of the source and target with the, summed
over all admissible {\em balanced} labellings of the 
internal edges with coefficents
equal to the product $\prod_{k=1}^{n-2} \Delta_{j_k}$, where the $n-2$
internal edges are labelled $(j_1,j_1), \ldots ,(j_{n-2},j_{n-2})$, and
the $n-1$ internal vertices are all evaluated a ``new 3-vertices''.

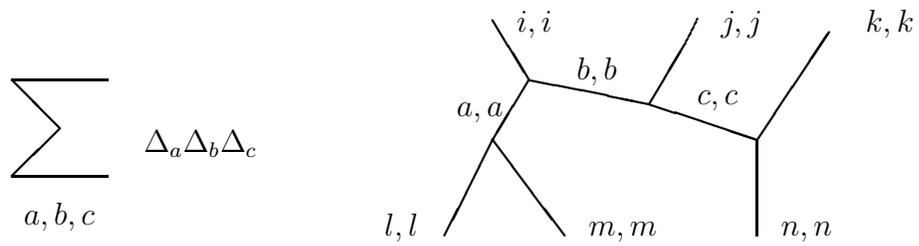
\begin{figure}[h]
\begin{center}

\setlength{\unitlength}{0.00083300in}%
\begingroup\makeatletter\ifx\SetFigFont\undefined
% extract first six characters in \fmtname
\def\x#1#2#3#4#5#6#7\relax{\def\x{#1#2#3#4#5#6}}%
\expandafter\x\fmtname xxxxxx\relax \def\y{splain}%
\ifx\x\y   % LaTeX or SliTeX?
\gdef\SetFigFont#1#2#3{%
  \ifnum #1<17\tiny\else \ifnum #1<20\small\else
  \ifnum #1<24\normalsize\else \ifnum #1<29\large\else
  \ifnum #1<34\Large\else \ifnum #1<41\LARGE\else
     \huge\fi\fi\fi\fi\fi\fi
  \csname #3\endcsname}%
\else
\gdef\SetFigFont#1#2#3{\begingroup
  \count@#1\relax \ifnum 25<\count@\count@25\fi
  \def\x{\endgroup\@setsize\SetFigFont{#2pt}}%
  \expandafter\x
    \csname \romannumeral\the\count@ pt\expandafter\endcsname
    \csname @\romannumeral\the\count@ pt\endcsname
  \csname #3\endcsname}%
\fi
\fi\endgroup
\begin{picture}(5337,1452)(1189,-1294)
\thicklines
\put(1801,-286){\line(-1, 0){600}}
\put(1201,-286){\line( 1,-1){300}}
\put(1501,-586){\line(-1,-1){300}}
\put(1201,-886){\line( 1, 0){600}}
\put(4201, 89){\line( 3,-5){225}}
\put(4426,-286){\line( 5,-1){750}}
\put(5176,-436){\line( 3, 5){225}}
\multiput(5401,-61)(3.75000,7.50000){21}{\makebox(6.6667,10.0000){\SetFigFont{7}{8.4}{rm}.}}
\put(4426,-286){\line(-3,-5){225}}
\put(4201,-661){\line( 3,-4){450}}
\put(4201,-661){\line(-1,-2){300}}
\put(3901,-1261){\line( 0, 1){  0}}
\put(5176,-436){\line( 3,-1){675}}
\put(5851,-661){\line( 2, 3){450}}
\put(5851,-661){\line( 0,-1){600}}
\put(4351, 14){\makebox(0,0)[lb]{\smash{\SetFigFont{12}{14.4}{rm}$i,i$}}}
\put(5626, 14){\makebox(0,0)[lb]{\smash{\SetFigFont{12}{14.4}{rm}$j,j$}}}
\put(6526, 14){\makebox(0,0)[lb]{\smash{\SetFigFont{12}{14.4}{rm}$k,k$}}}
\put(3526,-1261){\makebox(0,0)[lb]{\smash{\SetFigFont{12}{14.4}{rm}$l,l$}}}
\put(4801,-1261){\makebox(0,0)[lb]{\smash{\SetFigFont{12}{14.4}{rm}$m,m$}}}
\put(6001,-1261){\makebox(0,0)[lb]{\smash{\SetFigFont{12}{14.4}{rm}$n,n$}}}
\put(3976,-511){\makebox(0,0)[lb]{\smash{\SetFigFont{12}{14.4}{rm}$a,a$}}}
\put(4726,-286){\makebox(0,0)[lb]{\smash{\SetFigFont{12}{14.4}{rm}$b,b$}}}
\put(5476,-436){\makebox(0,0)[lb]{\smash{\SetFigFont{12}{14.4}{rm}$c,c$}}}
\put(2026,-736){\makebox(0,0)[lb]
{\smash{\SetFigFont{12}{14.4}{rm}$\Delta_a \Delta_b \Delta_c$}}}
\put(1276,-1186){\makebox(0,0)[lb]{\smash{\SetFigFont{12}{14.4}{rm}$a,b,c$}}}
\end{picture}

\caption{\label{a.6.vertex} A 6-vertex}
\end{center}
\end{figure}

An example is shown in Figure
\ref{a.6.vertex}.

It is then easy to apply Proposition \ref{turning.prop} to show

\begin{thm}
The value of an $n$-vertex depends only on the source and the target.
\end{thm}

\noindent{\bf proof:} Given any two trivalent trees with the same source
and target leaves, we can obtain one from the other by iteratively ``fusing''
edges as in Proposition \ref{turning.prop}. $\Box$
\smallskip

Once this is known, it is trivial to apply the observation that three-vertices
with balanced labels absorb braidings and are rotated by duality maps to 
justify the name of $n$-vertex for these operators.  In particular, recalling
that a prolongation of a map in a tensor category is an 
arbitary tensor product of
the map with identity maps, we have:

\begin{thm}
Any composition of a prolongation of an $n$-vertex with prolongations of
braidings, (associators, unit transformations,) and duality transformations
which admits a description as a connected and simply-connected 
recombination diagram when
the $n$-vertex is represented as a single $n$-valent vertex is itself
an $n$-vertex.
\end{thm}

\noindent{\bf sketch of proof:} Use the coherence theorem for tortile
categories \cite{shum} to rewrite all braidings and unit transformations
as compositions of prolongations of braidings, respectively unit 
transformations, indexed by single (balanced) irreducibles, not tensor
products.  For each braiding or duality transformation
composed directly with the $n$-vertex, use the previous proposition to
express the $n$-vertex in terms of a linear combination of trees so that
all objects indexing the braiding or unit transformation correspond
to arcs incident with a single trivalent vertex of the tree.  The braiding
or duality transformation may then be removed by applying the symmetry 
properties for 3-vertices noted above.$\Box$
\smallskip

What is not immediately clear is that there are also good canonical choices
for $0$-, $1$-, and $2$-vertices.  Observe that it is clear that a 
$0$ vertex must be an endomorphism of $(0,0)$, that is a scalar.  In 
considering possible relations with previously defined graph invariants,
it will be convenient to choose $1$.  For $1$- and $2$-vertices, however,
an additional structure present on the $n$-vertices already constructed
suggests a canonical choice in each case.

\begin{propo}
If $v_{x_1,\ldots ,x_k,r}^{y_1,\ldots ,y_{n-k-1}}$ denotes an $n$-vertex for
$n \geq 3$
with source $V_{(y_1,y_1)}\otimes \ldots \otimes V_{(y_{n-k-1},y_{n-k-1})}$
and target $V_{(x_1,x_1)}\otimes \ldots \otimes V_{(x_k,x_k)}\otimes 
V_{(r,r)}$, then

\[ \sum_r \Delta_r 
v_{x_1,\ldots ,x_k,r}^{y_1,\ldots ,y_{n-k-1}}(
V_{(x_1,x_1)}\otimes \ldots \otimes V_{(x_k,x_k)}\otimes v_{i,j}^r) \]

\noindent is an $n+1$-vertex.
\end{propo}

\noindent{\bf proof:} Immediate from the construction of $n$-vertices.
$\Box$
\smallskip

If we wish the result to hold without the restriction $n \geq 3$,
we must define $2$- and $1$-vertices by  

\[v_j^i = \left\{ \begin{array}{ll}\frac{1}{\Delta_j} Id_{V_{(j,j)}} & 
					\mbox{\rm if $i=j$} \\
				  0:V_{(i,i)}\rightarrow V_{(j,j)}  & 
					\mbox{otherwise}
		   \end{array} \right. \]

\noindent and

\[ v_j = \left\{ \begin{array}{ll}Id_{V_{(0,0)}} & 
					\mbox{\rm if $i=j$} \\
				  0:V_{(0,0)}\rightarrow V_{(j,j)}  & 
					\mbox{otherwise}
		   \end{array} \right. \]

\noindent respectively.

The check that these lead to the desired properties is easy and left to the
reader.

Thus we see that the full sub-category of tensor products of balanced 
irreducibles in $Rep(U_q^{bal}({\frak so}(4)))$ is a ``graphical category''
in the terminology of \cite{Y.graph}.
\newpage

\section{Invariants of Embedded Graphs}

Some preliminaries are in order:

\begin{defin}  A {\em graph} $\Gamma$ is a pair of sets $V(\Gamma), E(\Gamma)$
whose elements are called {\em vertices} and {\em edges} respectively, together
with a function 

\[ i_\Gamma:E(\Gamma) \rightarrow V(\Gamma)_1 \coprod V(\Gamma)_2. \]

\noindent where $S_i$ denotes the set of $i$-element subsets of $S$. 
Usually elements of the target is thought of as unordered pairs of elements
of $V(\Gamma)$.  By abuse of terminology and notation, we identify $\Gamma$ 
with its geometric realization, the topological space

\[  E(\Gamma)\times [0,1] / \equiv \]

\noindent where $\equiv$ is an equivalence relation induced by choosing
for each $e\in E(\Gamma)$ a surjection $k_e:\{0,1\}\rightarrow i_\Gamma(e)$, 
and defining $(e,x) \equiv (\eta, \xi)$ if and only if they are equal or
both $x$ and $\xi$ are $0$ or $1$ and $k_e(x) = k_\eta(\xi)$.
\end{defin}

Observe first that for any two choices of $\equiv$ there is a homeomorphism
between the resulting spaces which maps each $\{e\}\times [0,1]$ to itself
either by the identity map or by $t \mapsto 1-t$, and second that a graph 
has a natural PL-structure.

\begin{defin} An {\em embedded graph} is a PL-embedding of a graph into
${\Bbb R}^3$. Two embedded graphs are {\em equivalent} if they have isomorphic
underlying graphs and there exist a PL ambient isotopy of one to the other 
which preserves the set of vertices.
\end{defin}

It is fairly easy to see how to use the family of $n$-vertices constructed
above to construct invariants of embedded graphs:  choose a balanced
irreducible to label all the edges of the graph and a generic projection
of the graph, the invariant is obtained by evaluating the 
recombination network obtained by interpreting the graphical vertices
of valance $n$ as $n$-vertices in the algebraic sense, crossings as braiding,
and maxima and minima as duality transformations.  Isotopy invariance is
readily verified using the symmetry properties already established for 
$n$-vertices and the coherence theorem of Shum \cite{shum} for tortile
categories.

More generally, we may consider invariants of edge-colored graphs obtained
by choosing a balanced irreducible for each color and labelling all edges of
that color with the chosen irreducible.  In the following, however, we will
confine ourselves to invariants obtained by labelling all edges with the
same balanced irreducible.

\begin{defin}Let ${\cal G}_j[\Gamma](q)$ denote the value of
the invariant obtained by labelling all edges of an embedded graph
$\Gamma$ with the irreducible $(j,j)$.
\end{defin}

We can then easily show

\begin{propo} \label{odd.eulerian.only} 
For $j$ odd,   ${\cal G}_j[\Gamma](q) \neq 0$ implies $\Gamma$ is an
Eulerian graph, that is each component of $\Gamma$ is either a single 
vertex, or admits an Eulerian cycle (equivalently all valences are even).
\end{propo}

\noindent{\bf proof:} Any odd valence vertex will be mapped to an
algebraic $n$-vertex with all incoming labels $(j,j)$.  The underlying
Kauffman-Lins 3-vertices in any expression for the $n$-vertex have value
$0$. $\Box$
\smallskip

\begin{propo} If $\Gamma^\prime$ is obtained from $\Gamma$ by reversing
all crossings in a regular projection of $\Gamma$, then

\[ {\cal G}_j[\Gamma^\prime](q) = {\cal G}_j[\Gamma](q) \]
\end{propo}

\noindent{\bf proof:} Now, ${\cal G}_j[\Gamma](q)$ can be computed
by fusing arcs exiting each crossing (as in Kauffman and Lins   \cite{KL})
and summing over all (not necessarily
balanced) labellings of the  new arcs.  Using Kauffman and Lins  \cite{KL},
we may rewrite this sum as a linear combination of evaluations of 
planar networks at the cost of multiplying each summand by 
$A^{-k(k+2)+l(l+2)}$ (resp. $A^{k(k+2)-l(l+2)}$
whenever an arc labelled $(k,l)$ occurs at the
exit of a positive (resp. negative) crossing. Now observe that since the
value of a planar network is unchanged by replacing $q$ by $q^{-1}$, that
given any choice of labels for the new arcs, the value of the summand
in the expression for ${\cal G}_j[\Gamma]$ corresponding to this choice
is equal to the value of the summand in the expression for 
${\cal G}_j[\Gamma^\prime]$ in which the first and second coordinates of
each label have been interchanged.  Since interchanging coordinates of 
labels is an involution on the the set of labellings of the new arcs,
the desired result follows. $\Box$

\smallskip  

 Now, letting $J[L](q)$ denote the Jones polynomial of a link $L$ we have

\begin{propo} If $\Gamma$ is a disjoint union of cycles, let $L$ be the
underlying link obtained by forgetting the vertices, equipped with any
orientation, and $v_\Gamma$ be the number of vertices of $\Gamma$, then

\[ {\cal G}_1[\Gamma](q) = \frac{J[L](q)J[L](q^{-1})}{\Delta_1^{v_\Gamma}} \]
\end{propo}

\noindent{\bf proof:} The value of the (1,1)-labelled network is 
$J_{f}[L](q)J_{f}[L](q^{-1})$ where $J_f$ denotes the framed Jones 
polynomial (Kauffman bracket, but written in $q$ rather than $A$). But the
normalization constants cancel, so this is the desired numerator.  The
definition of the 2-vertex contributes a factor of 
$\frac{1}{\Delta_1^{v_\Gamma}}$. $\Box$

\smallskip

	A standard notion in graph theory is that of a {\em cutpoint}, that is
a vertex whose removal increases the number of connected
components. This notion is usually defined in a purely combinatorial way so
that the presence of loops at a vertex does not imply that the vertex is a
cutpoint.  We will need a stronger notion of of cutpoint which takes into
account the topology and embedding of the graph.  First let us recall
a standard notion from link theory (in a form applicable to embedded graphs):

\begin{defin}
An embedded graph $\Gamma$ is a {\em separated union} if it is ambient isotopic
to an embedded graph $\eta:\Gamma_1 \coprod \Gamma_2 \rightarrow {\Bbb R}^3$
where the images of $\Gamma_1$ and $\Gamma_2$ lie in disjoint open balls.
\end{defin}

Similarly we make

\begin{defin} An embedded graph $\Gamma$ is an {\em almost-separated union} of
two subgraphs $\Gamma_1$ and $\Gamma_2$ if the intersection of $\Gamma_1$ and
$\Gamma_2$ is a single vertex $v$, and $\Gamma$ is ambient isotopic to
a graph in which the images of $\Gamma_1\setminus\{v\}$ and 
$\Gamma_2\setminus\{v\}$ lie in disjoint open balls.
\end{defin}

Note that here we must regard graphs as their geometric realization, since by
$\Gamma_i\setminus\{v\}$ we mean the space obtained by deleting the point
$v$, not the graph obtained by deleting $v$ and all incident edges.

\begin{defin}
A {\em topological cutpoint} of an embedded graph $\Gamma$ is a vertex $v$ 
such that $\Gamma$ is an almost-separated union of two subgraphs $\Gamma_1$ and
$\Gamma_2$ with common
vertex $v$.  By a {\em splitting} of a graph at a topological cut-point, we
mean the separated union $\Gamma^\prime$ of $\Gamma_1$ and $\Gamma_2$.
\end{defin}

\begin{thm}
If $v$ is a topological cutpoint of $\Gamma$ and $\Gamma^\prime \cong
\Gamma_1 \coprod \Gamma_2$ is a splitting of $\Gamma$ at $v$, then

 \[{\cal G}[\Gamma] = {\cal G}[\Gamma^\prime] = 
{\cal G}[\Gamma_1] {\cal G}[\Gamma_2] \]
\end{thm}

\noindent{\bf proof: } Observe that is $v$ is a topological cutpoint, we may
construct the tree over which the algebraic vertex assigned to $v$ is defined
in such a way that there is a distinguished edge whose removal 
separates $\Gamma_1$ and
$\Gamma_2$. Further observe that if we use this tree to describe the 
algebraic vertex at $v$, and render the evaluation of ${\cal G}[\Gamma]$
as a sum of recombination diagrams, all of the diagrams have turn-arounds
incident with distinguished edge.  By one of the properties of 
$q$-symmetrizers noted in \cite{KL}, it follows that only those summands
in which the label on the distinguished edge is $(0,0)$ can be non-zero.

It may readily be verified that the value of the summand in a 
algebraic vertex with
label $(0,0)$ on one edge in the tree describing the algebraic vertex
is a tensor product of two algebaric vertices whose trees are those
obtained by removing the distinguished edge.  The result then follows
from the lemma below. $\Box$  

\begin{lemma}
${\cal G}_j$ is multiplicative under separated union.
\end{lemma}

\noindent{\bf proof:} Immediate from the functorial construction.$\Box$

\newpage
\section{Conclusions}

The preceding lemma and Proposition \ref{odd.eulerian.only} are strongly 
suggestive that (for odd $j$) ${\cal G}_j$ may well be related to the
Martin polynomial (cf. \cite{M}, \cite{E-M}), at least at 
$q = \pm 1$ where the braiding is trivial.  It would seem to be 
a fruitful undertaking to examine the invariants described herein in
terms of the Hopf algebra structure on the space of linear combinations
of embedded graphs induced by the formulas given in \cite{E-M} for 
(non-embedded) graphs.


\newpage

\begin{thebibliography}{99}

\bibitem{BC}
Barrett, J.W. and Crane, L., ``Relativistic Spin Networks and Quantum 
Gravity,'' gr-qc 9709028 (to appear CQG).

\bibitem{CKY}
Crane, L., Kauffman, L.H., and Yetter, D.N., ``State-Sum Invariants of
4-Manifolds,'' {\em JKTR} {\bf 4} (2) (1997) 177-234.

\bibitem{E-M}
Ellis-Monaghan, J.A., ``New Results for the Martin Polynomial,'' preprint
(1997). 

\bibitem{FY}
Freyd, P.J. and Yetter, D.N., ``Braided Compact Closed Categories with
Applications to Low-Dimensional Topology,''  {\em Adv. in Math.} {\bf 77}
(2) (1989) 156-182.

\bibitem{JS}
Joyal, A. and Street, R., ``The Geometry of Tensor Calculus, I,'' 
{\em Adv. in Math.} {\bf 88} (1) (1991) 55-112.

\bibitem{KL}
Kauffman, L.H. and Lins, S.L., {\em Temperley-Lieb Recoupling Theory
and Invariants of 3-Manifolds,} Princeton Univ. Press (1994).

\bibitem{M}
Martin, P. {\em Enumerations Euleriennes dans les Multigraphs et 
Invariants de Tutte-Grothendeick}, Thesis, Grenoble, 1977. 

\bibitem{CWM} Mac Lane, S. {\em Categories for the Working Mathematician},
Springer (1971).

\bibitem{RP} 
G. Ponzano and T. Regge, ``Semiclassical Limit of Racah
Coefficients'' in {\em Spectroscopic and Group Theoretic Methods in Physics}
(R. Bloch, ed), North Holland, Amsterdam, 1968.

\bibitem{shum} Shum, M.-C. ``Tortile Tensor Categories,'' {\em J. Pure and
App. Alg.} {\bf 93} (1) (1994), 57-110.

\bibitem{Y.rep}
Yetter, D.N., ``Quantum Groups and Representations of Monoidal Categories,''
{\em Math. Proc. Camb. Phil. Soc.} {\bf 108} (1990) 261-290.

\bibitem{Y.graph}
Yetter, D.N., ``Category Theoretic Representations of Embedded Graphs in
${\bf S}^3$,'' {\em Adv. in Math} {\bf 77} (2) (1989) 137-155. 


\end{thebibliography}
\end{document}